\documentclass{article}
\usepackage{amssymb,amsmath}
\usepackage{graphicx}
\usepackage{ dsfont }

\def\qed{\hfill {\hbox{${\vcenter{\vbox{               
   \hrule height 0.4pt\hbox{\vrule width 0.4pt height 6pt
   \kern5pt\vrule width 0.4pt}\hrule height 0.4pt}}}$}}}

\def\hat{\widehat}

\newtheorem{theorem}{Theorem}
\newtheorem{definition}{Definition}

\newtheorem{corollary}[theorem]{Corollary}
\newtheorem{example}{Example}
\newtheorem{remark}[example]{Remark}

{\qed\par\medskip}

\date{}

\title{\Large \textbf{Polynomial knot and link invariants from the virtual 
biquandle}}

\author{Alissa S. Crans\footnote{acrans@lmu.edu, Department of Mathematics, Loyola Marymount University, One LMU Drive, Suite 2700, Los Angeles, CA 90045}
\and Allison Henrich \footnote{henricha@seattleu.edu,
Mathematics Department, Seattle University, 901 12th Avenue, PO Box 222000, Seattle, WA 98122} 
\and Sam Nelson\footnote{knots@esotericka.org, Department of Mathematical Sciences, Claremont McKenna College, 850 Columbia Ave., Claremont, CA 91711
}}

\begin{document}
\maketitle

\begin{abstract}
The Alexander biquandle of a virtual knot or link is a
module over a 2-variable Laurent polynomial ring which is an 
invariant of virtual knots and links. The elementary ideals of this 
module are then invariants of virtual isotopy which determine 
both the generalized Alexander polynomial (also known as the Sawollek 
polynomial) for virtual knots and the classical Alexander polynomial 
for classical knots. For a fixed monomial ordering $<$, the Gr\"obner 
bases for these ideals are computable, comparable invariants which fully
determine the elementary ideals and which generalize and unify
the classical and generalized Alexander polynomials. We provide examples 
to illustrate the usefulness of these invariants and propose questions 
for future work.
\end{abstract}

\parbox{5.5in} {\textsc{Keywords:} virtual knot, generalized Alexander 
polynomial, virtual Alexander polynomial, Sawollek polynomial, biquandle, 
Alexander biquandle.
\smallskip

}

\section{\large \textbf{Introduction}}

The Alexander biquandle of an oriented classical or virtual knot
or link is a module over a ring of 2-variable Laurent polynomials which is 
invariant under virtual Reidemeister moves. In this paper we describe
methods of obtaining computable invariants of classical and virtual knots
and links from the Alexander biquandle. These families of invariants
generalize and unify the classical Alexander polynomials for classical knots
and the 
generalized Alexander polynomial for virtual knots and links (also known
as the \textit{virtual Alexander polynomial} \cite{K, KR, KM}  and the 
\textit{Sawollek polynomial} \cite{S}), $Z_{D}(x,y)$.

This work was inspired by the Remark following Theorem 3 in Sawollek's 
paper \cite{S} which 
states, ``For a connected sum $D_1 \# D_2$ of virtual link diagrams $D_1$ 
and $D_2$, a formula of the form
\[Z_{D_1 \# D_2} (x,y) = c Z_{D_1} (x,y) Z_{D_2}(x,y)\]
with a constant $c$ does not hold in general..."
This is not surprising since the connected sum of virtual knots and 
links is not well-defined, yet in certain cases the result nevertheless holds.  
Computing the virtual Alexander polynomial of various connected sums of 
virtual knots reveals that, as anticipated, these polynomials 
differed depending on where the connect sum is performed.  
Moreover, these 
polynomials gave no indication of where the connect sum was taken. It is 
natural, then, to ask under what circumstances does 
$Z_{D}(x,y)$ satisfy the above equation.

Since the classical Alexander polynomial
generates the $k=1$ elementary ideal of the Alexander matrix, we expect 
the factor relationship of the polynomials of connected summands to result in 
an inclusion
relationship at the level of ideals. To begin our study of such relationships,
therefore, we turned our attention to describing these ideal-valued invariants 
in more detail as a first step. The result is the present paper; we anticipate 
future work taking these ideas further.

The paper is organized as follows. In Section \ref{VAQ} we recall the 
Alexander biquandle. In Section \ref{P} we define invariants of oriented
virtual and classical knots and links derived from the Alexander  
biquandle using principal ideals. In Section \ref{G} we use reduced 
Gr\"{o}bner bases to define further invariants from the Alexander 
biquandle. In Section \ref{CE} we collect some computations and examples. 
We conclude with some questions for future work in Section \ref{Q}.

\section{\large \textbf{The Alexander biquandle}}\label{VAQ}

Recall from \cite{K} that an \textit{oriented virtual link} is an equivalence
class of \textit{oriented virtual link diagrams} under the equivalence relation
defined by the \textit{virtual Reidemeister moves}, obtained by considering
all possible oriented versions of the moves pictured below:
\[\includegraphics{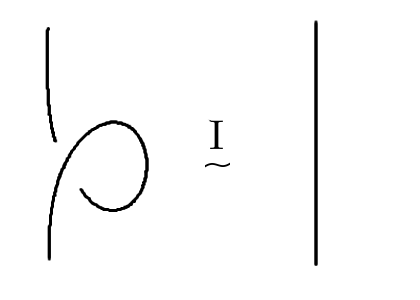} \quad
\includegraphics{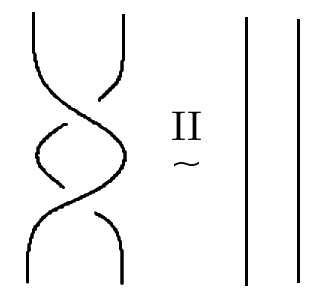} \quad
\includegraphics{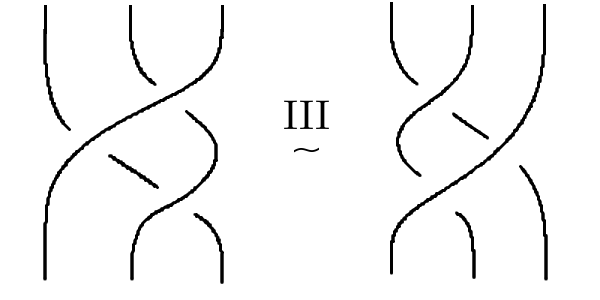}
\]
\[
\includegraphics{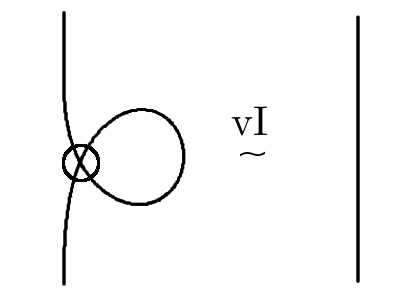} \quad
\includegraphics{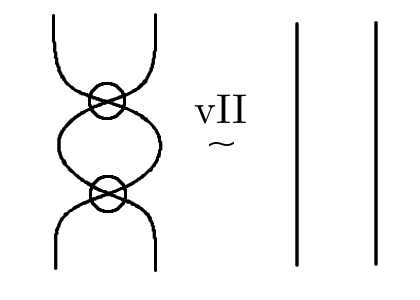} \]\[
\includegraphics{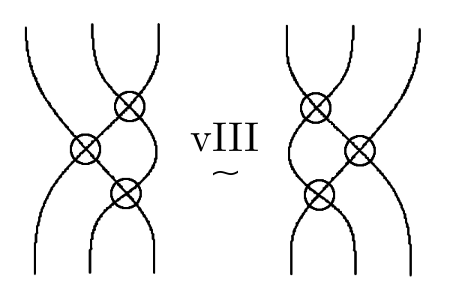} \quad
\includegraphics{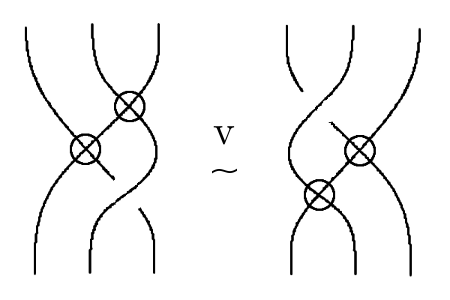}
\]

We now review the definitions of the virtual biquandle and the fundamental 
virtual biquandle. This second definition will allow us to see the 
relationship between biquandles and virtual knots.

\begin{definition}\textup{Let $X$ be a set and define the \textit{diagonal map} 
$\Delta:X\to X\times X$
by $\Delta(x)=(x,x)$. A \textit{virtual biquandle}
structure on $X$ consists of two invertible maps 
$B,V:X\times X\to X\times X$ satisfying the axioms:
\begin{itemize}
\item[(0)] $V^2=\mathrm{Id}:X\times X\to X\times X$,
\item[(1)] There exist unique invertible maps $S,vS:X\times X\to X\times X$ 
(called the \textit{sideways maps}) satisfying
\[S(B_1(x,y),x)=(B_2(x,y),y) \quad \mathrm{and} \quad
vS(V_1(x,y),x)=(V_2(x,y),y)\]
for all $x,y\in X$,
\item[(2)] 
$(S^{\pm 1}\circ \Delta)_{j}$ and $(vS^{\pm 1}\circ \Delta)_{j}$ are bijections
for $j=1,2$ satisfying\[
(S\circ \Delta)_1=(S\circ \Delta)_2 
\quad \mathrm{and} \quad
(vS\circ \Delta)_1=(vS\circ \Delta)_2, \; \textrm{and} \] 
\item[(3)] $B$ and $V$ satisfy the \textit{set-theoretic Yang-Baxter 
equations}
\[\begin{array}{rcl}
(B\times I)(I\times B)(B\times I) & = & (B\times I)(B\times I)(I\times B) \\
(V\times I)(I\times V)(V\times I) & = & (V\times I)(V\times I)(I\times V) \\
(V\times I)(I\times B)(V\times I) & = & (V\times I)(B\times I)(I\times V) \\
\end{array}\]
\end{itemize}
}\end{definition}

We note that these axioms result from by applying the labeling condition 
pictured below to the virtual Reidemeister moves. For 
instance, the invertibility of $B$ and $V$ along with Axiom (1) represent 
Reidemeister 2 moves. Reidemeister 1 moves are represented by Axiom (2), and 
Reidemeister 3 moves are related to Axiom (3). 
%

Suppose we are given a virtual knot or link $L$ with fixed diagram. Let $X=\{x_1,\dots,x_n\}$
be a set of generators that are in one-to-one correspondence with the set of 
\textit{semiarcs} in $L$, i.e. the portions of the virtual knot or link between
crossing points (whether virtual, classical over or classical under). The set 
of \textit{virtual biquandle words in $X$}, denoted $W(X)$, is defined 
recursively by the rules:
\begin{itemize}
\item $x\in X$ implies $x\in W(X)$ and
\item $x,y\in W(X)$ implies $B^{\pm 1}_{j}(x,y),\ V^{\pm 1}_{j}(x,y),\ 
S^{\pm 1}_{j}(x,y),\; \textrm{and} \; vS^{\pm 1}_{j}(x,y)\in W(X)$ for $j=1,2$. 
\end{itemize}

The \textit{free virtual biquandle on $X$}, denoted $FV(X)$, is the set of
equivalence classes of virtual biquandle words in $X$ under the equivalence 
relation on $W(X)$ generated by the virtual biquandle axioms. 

We remark that $FV(X)$ gives provides almost no information about our
original knot or link $L$. In order to capture information contained in 
our fixed diagram of $L$, we mod out by relations suggested by the crossings
to obtain the fundamental virtual biquandle.

The \textit{fundamental virtual biquandle of $L$}, denoted $FB(L)$, is the 
set of equivalence classes of $FV(X)$ under the equivalence relation 
generated by the \textit{crossing relations} in $L$:
\[\begin{array}{cc}
\includegraphics{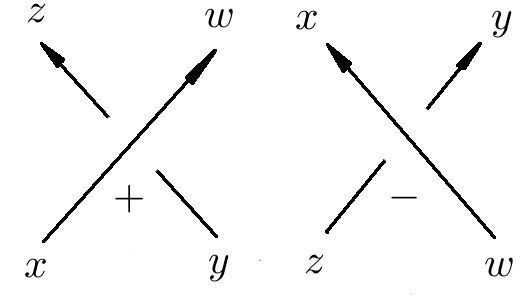}  & \qquad \qquad
\includegraphics{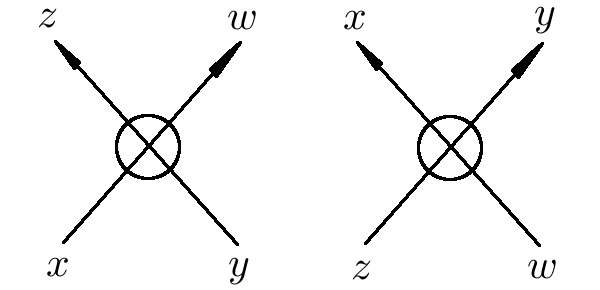}\\
B(x,y) = (z,w) &\qquad \qquad V(x,y)=(z,w)
\end{array}
\]
Alternatively, we can 
describe $FB(L)$ more directly as the set of equivalence classes 
in $W(X)$ under the equivalence relation generated by both the crossing 
relations and the virtual biquandle axioms.

\begin{theorem} Let $L$ be a virtual link written as a closed virtual braid 
$\beta$, that is
$L=\hat{\beta}$. The fundamental virtual biquandle of a virtual link $L$ is 
isomorphic to the fundamental virtual biquandle of the link $L'$ obtained 
from $L$ by 
reversing the direction of all strands in the closure of the inverse braid
$\hat{\beta^{-1}}$.
\end{theorem}

\noindent The proof is effectively the same as the analogous result for 
classical 
biquandles of virtual links in \cite{KH}; we provide an illustration.
Consider the virtual link $L_1=\hat{\beta}$ given by the closure of the braid 
\[\raisebox{-0.65in}{\includegraphics{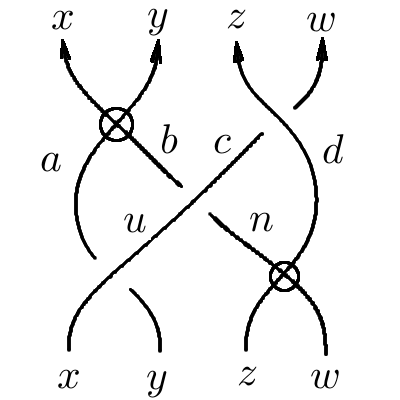}}
 FB(L_1) = \left\langle \begin{array}{l|l} 
x,y,z,w,u, & V(a,b)=(x,y), \; \; B(z,w)=(c,d), \\ 
n,a,b,c,d &  B(u,n)=(b,c), \; \; B(x,y)=(a,u), \\ &  V(z,w)=(n,d) 
 \\\end{array}  \right\rangle
\]
The reversed inverse braid closure $L_2=r\hat{\beta^{-1}}$ has isomorphic 
virtual biquandle (indeed, with identical presentation) to that of $L_1$ 

\[\raisebox{-0.65in}{\includegraphics{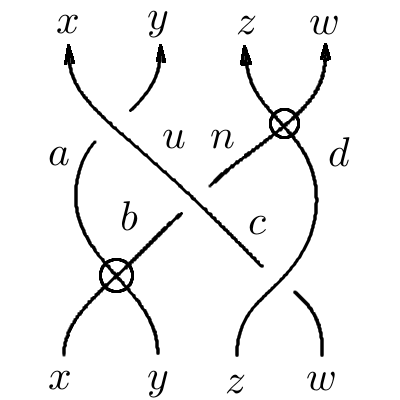}}
 FB(L_2) = \left\langle \begin{array}{l|l}
x,y,z,w,u, & V(a,b)=(x,y), \; \; B(z,w)=(c,d), \\ 
n,a,b,c,d &  B(u,n)=(b,c), \; \; B(x,y)=(a,u), \\ &  V(z,w)=(n,d)
 \\ \end{array} \right\rangle \]
while the unreversed inverse braid closure $L_3=\hat{\beta^{-1}}$ has a 
generally distinct virtual biquandle:

\[\raisebox{-0.65in}{\includegraphics{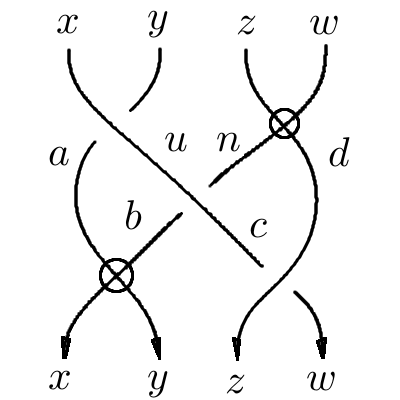}}
FB(L_3) = \left\langle \begin{array}{l|l}
x,y,z,w,u, & B(u,a)=(y,x), \; \; V(w,z)=(d,n), \\ 
n,a,b,c,d  & B(n,u)=(c,b), \; \; V(b,a)=(y,x), \\ &  B(w,z)=(d,c) 
 \\ \end{array} \right\rangle \] 
 
\begin{remark}\textup{
The reversed inverse braid closure $r\hat{\beta^{-1}}$ is called the 
\textit{vertical mirror image}
in \cite{KH} and is one of the $2^c$ possible orientations for a $c$-component 
link of what the Knot Atlas \cite{G} calls the \textit{horizontal mirror image}.
Another ``mirror image'' of a virtual link can be obtained by switching every
overcrossing to an undercrossing while fixing the diagram outside a 
neighborhood of the crossing. This operation is called the \textit{mirror 
image} in \cite{KH} and is one orientation of the Knot Atlas' \textit{vertical 
mirror image}. To avoid confusion we will refer to this operation as the 
\textit{sign switch} of $L$, denoted $sL$; it is in general a distinct 
oriented virtual link from $L:$}
\[\raisebox{-0.65in}{\includegraphics{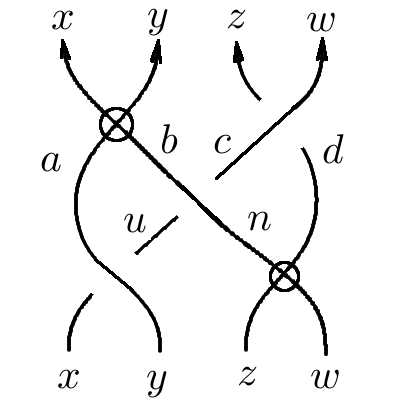}}
 FB(sL) = \left\langle  \begin{array}{l|l}
x,y,z,w,u,  & B(a,u)=(x,y), \; \; V(z,w)=(n,d), \\ 
n,a,b,c,d &   B(b,c)=(u,n), \; \; V(a,b)=(x,y), \\ &  B(c,d)=(z,w) 
 \\\end{array} \right\rangle\]
\end{remark}

An important example of a virtual biquandle, and one to which we will devote 
considerable attention, is the Alexander biquandle. We begin with 
some notation. Let $\Lambda=\mathbb{Z}[t^{\pm 1},s^{\pm 1}]$ be the
ring of Laurent polynomials in two variables $t,s$. We can 
think of  $\Lambda$ as a quotient ring of the four-variable polynomial ring
$\tilde\Lambda=\mathbb{Z}[t,s,t^{-1},s^{-1}]$ by the 
ideal $I=\langle tt^{-1}-1,ss^{-1}-1\rangle$

\begin{definition}\textup{
Let $L$ be an oriented virtual knot or link diagram and let 
$X=\{x_1,\dots,x_n\}$ be a set of generators corresponding to the 
semiarcs of $L$. Then the
\textit{Alexander biquandle} of $L$, denoted $AB(L)$, is the 
$\Lambda$-module generated by $X$ with the relations pictured below at 
positively ($+$) and negatively ($-$) oriented classical crossings and 
at virtual crossings:
\[\begin{array}{cc}
\quad \quad \includegraphics{ac-ah-sn-1.png} \quad \quad  & \quad  \quad
\includegraphics{ac-ah-sn-2.png} \quad\quad\\
\begin{array}{rcl}
z & = & ty+(1-st)x \\
w & = & sx
\end{array} & 
\begin{array}{rcl}
z & = & y \\
w & = & x
\end{array} 
\end{array}
\]
In particular, $AB(L)$ is obtained from the fundamental virtual biquandle
of $L$ by setting \[B(x,y)=(ty+(1-st)x,sx)\quad \mathrm{and}\quad
V(x,y)=(y,x).\]
}\end{definition}

\noindent A straightforward check verifies that the Alexander biquandle of a 
virtual knot is invariant under virtual Reidemeister moves; see \cite{KM}
for details.

\begin{remark}
\textup{It is natural to consider, as we initially did, including a 
coefficient at the virtual crossings, e.g. set $V(x,y)=(vy,v^{-1}x)$. 
However, Theorem 7.1 in \cite{B} implies that the resulting 
$\mathbb{Z}[t^{\pm 1},s^{\pm 1}, v^{\pm 1}]$-module contains the same 
information as the simpler $\mathbb{Z}[t^{\pm 1},s^{\pm 1}]$-module. The 
authors would like to thank Lou Kauffman for bringing this result to 
their attention.}
\end{remark}

We can represent the Alexander biquandle of a virtual knot or link 
with the coefficient matrix of the homogeneous system of equations 
determined by the crossings, known as a \textit{presentation matrix}.

\begin{example}\label{ex1}\textup{
The virtual knot below has Alexander biquandle with presentation and
presentation matrix as listed.}
\[\raisebox{-0.8in}{\scalebox{1.3}{\includegraphics{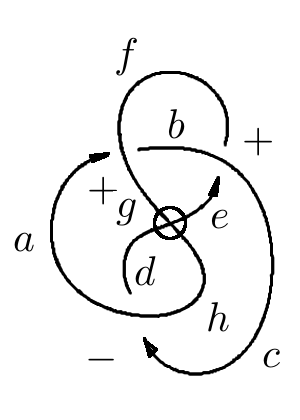}}}\quad\quad
\begin{array}{c}
AB(L)= \left\langle 
\begin{array}{r|l}
a,b,c,d,e,f,g,h &
\begin{array}{rcl}
 b & = & ta+(1-st)f, \\
 g & = & sf, \\ 
 f & = & te+(1-st)b, \\
 c & = & sb, \\ 
 e & = & d, \\
 g & = & h, \\ 
 c & = & td+(1-st)a,\\ 
 h & = & sa \end{array}
\end{array} \right\rangle\\
 \\
M=\left[\begin{array}{cccccccc}
t & -1 & 0 & 0 & 0 & 1-st & 0 & 0 \\
0 &  0 & 0 & 0 & 0 & s & -1 & 0 \\
0 & 1-st & 0 & 0 & t & -1 & 0 & 0 \\
0 & s & -1 & 0 & 0 & 0 & 0 & 0 \\
0 & 0 & 0 & 1 & -1 & 0 & 0 & 0 \\
0 & 0 & 0 & 0 & 0 & 0 & -1 & 1 \\
1-st & 0 & -1 & t & 0 & 0 & 0 & 0 \\
s & 0 & 0 & 0 & 0 & 0 & 0 & -1 
\end{array}\right]\\
\end{array}\]
\textup{For comparison, the Sawollek \cite{S} matrix $M-P$ is:}
\[
 M - P = \left[ \begin{array}{cccccc}
1 - x & -y & 0 & 0& -1 & 0 \\
-x/y & 0 & 0 & -1 & 0 & 0 \\
0 & 0 & 1-x & -y & 0 & -1 \\
0 & -1 & -x/y & 0 & 0 & 0 \\
0 & 0 & -1 & 0 & 0 & -y/x \\
-1 & 0 & 0 & 0 & -1/y & 1- 1/x 
\end{array} \right]. \]
\textup{We describe the relationship between these matrices in 
Example~\ref{corresp}.
}
\end{example}

\begin{definition}\textup{
Let $L$ be a virtual link, $AB(L)$ the Alexander biquandle of
$L$, and $M$ an $m\times n$ presentation matrix for $AB(L)$. 
The ideal $I_k$ of 
$\Lambda$ generated by the $(m-k)$-minors of $M$ is the 
\textit{$k$th elementary ideal} of $M$. 
}\end{definition}

It is well-known that:
\begin{theorem}
The elementary ideals 
of a module over a commutative ring
with identity do not depend on the choice of presentation matrix 
for the module.
\end{theorem}

\noindent In particular, any two presentation matrices of the same module 
differ by a sequence of moves of the forms:
\begin{itemize}
\item reordering of rows or columns,
\item adding or deleting a row of all zeroes,
\item adding or deleting a row and column with a $1$ in the intersection and
all other entries $0$,
\item adding a scalar multiple of one row (or column, respectively) to
another row (or column, respectively), or
\item replacing a row or column by an invertible scalar multiple of itself.
\end{itemize}

\noindent 
One then checks that these moves do not change the ideal. For instance, 
switching the order of two rows will multiply the minors by $-1$, but
this does not change the ideal since $-1$ is a unit. Similarly, 
multilinearity of the determinant ensures that adding a scalar multiple of
one row to another does not change the minors, etc. See \cite{C,L} for more
details.

It then follows that:
\begin{corollary}
Let $L$ be a virtual link and let $M$ be a presentation matrix for
$AB(L)$. Then the elementary ideals $I_k$ of $M$ are invariants of virtual 
isotopy.
\end{corollary}

In order to compare ideals, we need some machinery. In the next sections we
describe two methods for comparing the ideals $I_k$.

\section{\large \textbf{Principal Alexander Polynomials}}\label{P}

In this section we employ the tried and true method of obtaining invariants
from $I_k$ using principal ideals, i.e. ideals generated by a single
element of $\Lambda = \mathbb{Z}[t^{\pm 1},s^{\pm 1}]$. 

\begin{definition}
\textup{Let $L$ be a virtual link. For any $k=0,1,2,\dots$, the 
\textit{$k$th principal Alexander polynomial} of $L$, denoted 
$\Delta_k^p(L)$, is the generator of the smallest principal ideal 
$P_k\subset \Lambda$ containing the $k$th elementary ideal $I_k$ of $AB(L)$.
}\end{definition}

\begin{remark}\textup{
Note that $\Delta_k^p(L)$ is only defined up to multiplication by units in
$\Lambda$, so two values $f$ and $g$ of $\Delta_k^p(L)$ are equivalent if 
$f=\pm t^ns^m g$ for some $n,m\in \mathbb{Z}$.
}\end{remark}

The classical Alexander polynomials can be obtained by the analogous 
construction starting with a presentation matrix for the the Alexander 
quandle, which is the special case of the Alexander biquandle with 
\[B(x,y)=(ty+(1-t)x,y)\quad\mathrm{and}\quad V(x,y)=(y,x)\]
or equivalently the result of specializing $s=1$ in $\Delta_k^p(L)$. 
In the literature, $\Delta_0^p$ is known as the \textit{generalized 
Alexander polynomial} or \textit{virtual Alexander polynomial}, and,
after a change of variables, the \textit{Sawollek polynomial} 
\cite{K,KR,KM}. In particular, specializing $s=1$ in $\Delta_1^p(L)$ for 
classical links yields the classical Alexander polynomial $\Delta(L)$.

\begin{example}\label{corresp}
\textup{The virtual knot in Example \ref{ex1} has 0th principal polynomial}
\[\Delta_0^p(K)=(1-st)(1-s)(1-t)\]
\textup{Taking the determinant of the Sawollek matrix in Example \ref{ex1}, 
performing the change of variables $x = st$ and $y = -s$ and canceling units
gives the generalized Alexander polynomial in the notation of \cite{K}:}
\[G_{K}(s,t)=(s - 1)(t - 1)(st - 1)=\Delta_0^p(K).\]
\textup{The classical Alexander polynomial of this virtual knot is 
$\Delta(K)=1.$}
\end{example}

Since a virtual knot always has $AB(L)$ presented by a square matrix, 
the top level $k=0$ ideal is always principal, generated by the determinant 
of the presentation matrix. To compute $\Delta^p_k(L)$ for $k>0$, we can 
find the $(m-k)$-minors of the presentation matrix, and after multiplying
by units if necessary to get elements of the UFD (Unique Factorization Domain) $\mathbb{Z}[t,s]$, find
the greatest common divisor.

\section{\large \textbf{Alexander-Gr\"{o}bner Invariants}}\label{G}

The principal polynomials $\Delta^p_k$ have the advantages 
of being fairly quick to compute and easy to compare; however, for any
principal ideal $P\subset\Lambda$, there are potentially many distinct 
non-principal 
ideals contained in $P$, and thus passing from $I_k$ to $P_k$ represents
a loss of information. To avoid this loss of information, we can employ the
idea of a Gr\"{o}bner basis.

Briefly, in a multivariable polynomial ring $R[x_1,\dots,x_n]$ over a PID (Principal Ideal Domain) $R$,
a \textit{term ordering} is a well-ordering on the set of monomials; 
such an ordering then gives each polynomial a well-defined 
\textit{leading term.} Standard examples of term orderings include:
\begin{itemize}
\item \textit{lexicographical ordering.} Starting with an ordering on the
variables, e.g. $x_1<x_2<\dots<x_n$, one compares terms by comparing powers
on the variables in order, with ties in one variable resolved by comparing 
the next. For example, if $x<y<z$ then we have
\[xy^2z<xy^2z^2<xy^3<x^2.\]
\item \textit{graded lexicographical ordering.} Here we start by comparing
total degree, with ties broken lexicographically. Thus for the previous
example in graded lexicographical ordering we have
\[x^2<xy^2<xy^3<xy^2z^2.\]
\end{itemize}

A generating set $G$ for an ideal $I$ in a polynomial ring is a 
\textit{Gr\"obner basis} for $I$ with respect to a choice of term ordering $<$
if the leading term of every element of $I$ is divisible by the leading term
of some element of $G$. Gr\"obner basis have many useful properties, such as:
\begin{itemize}
\item reducing a polynomial $f\in I$ by $G$ using the multivariable division 
algorithm (i.e. repeatedly subtracting from $f$ multiples of elements of
$G$ to cancel leading terms) always results in $0$ (this is not true for 
arbitrary generating sets), and
\item the remainder of a polynomial after multivariable division by $G$ is 
unique, enabling computations in the quotient ring $\Lambda/I$.
\end{itemize}
Moreover, a Gr\"obner basis is \textit{reduced} if redundant elements have been
removed, i.e. if every leading term has coefficient $1$ and no monomial in any
element $g\in G$ is in the ideal generated by the leading terms of 
$G - \{g\}$. Importantly for us, for a given choice of term ordering, 
reduced Gr\"obner bases of ideals are unique.

Thus, we would like to compare the ideals $I_k$ by finding and comparing 
Gr\"obner bases. One slight problem is that 
$\Lambda = \mathbb{Z}[t^{\pm 1},s^{\pm 1}]$ 
is a ring of \textit{Laurent} polynomials, not polynomials. To address this, 
we will pull 
back from the ring $\Lambda$ to the four-variable polynomial ring 
$\tilde\Lambda=\mathbb{Z}[t,s,t^{-1},s^{-1}]$ in which $t^{-1}$ and $s^{-1}$
are considered new variables as opposed to powers of  $t$ and $s$.

\begin{definition}
\textup{Let $L$ be a virtual link and $M$ a presentation matrix for 
$AB(L)$. For each elementary ideal $I_k$, let $\tilde{I_k}$ be the preimage
of $I_k$ in $\tilde\Lambda$. For a choice of term ordering $<$, 
the $k$th \textit{Alexander-Gr\"obner Invariant} $\Delta^<_k(L)$ is the 
Gr\"obner basis of $\tilde\Lambda$ with respect to the term ordering $<$.}
\end{definition}

To compute $\Delta^<_k(L)$, one uses \textit{Buchberger's algorithm} \cite{E} starting
by setting $G$ equal to the set of $(m-k)$-minors together with the 
generators $1-tt^{-1}$ and $1-ss^{-1}$. For each pair of elements 
$g,g'\in G$, we find a difference $S=pg-qg'$ using the least common multiple 
of the leading terms of $g$ and $g'$ such that that the leading terms of 
$pg$ and $qg'$ cancel. After reducing $S$ mod $G$,
if the remainder is nonzero, add it to $G$ and start over. When no more
polynomials are added, we have a Gr\"obner basis. We then remove elements
whose monomials lie in the ideal generated by the leading terms of the
other elements of $G$ to obtain the reduced Gr\"obner basis, unique up
to choice of term ordering $<$. See \cite{E} for more.

Since the resulting sets of polynomials can be quite large, we may employ
various strategies to obtain more conveniently comparable invariants at the 
cost of losing some information. These include but are certainly not limited to:
\begin{itemize}
\item \textit{Alexander-Gr\"obner cardinality.}  The number of elements in 
the Gr\"obner basis of $\tilde{I}_k$, 
$\Delta^{<,\mathbb{Z}}_k(L)=|\Delta^<_k(L)|$;
\item \textit{Alexander-Gr\"obner sum.} The sum of the elements of 
$\Delta^<_k(L)$, \[\Delta^{<,\Sigma}_k(L)=\sum_{g\in \Delta^<_k} g;\] 
and
\item \textit{Alexander-Gr\"obner maximum polynomial.} The maximal
element of $\Delta^<_k(L)$ with respect to the term ordering $<$, 
$\Delta^{<,M}_k(L)$.
\end{itemize}

We note that the last two invariants coincide with the principal polynomials
when the ideals $I_k$ are principal; while the first invariant is $1$ iff the
ideal $I_k$ is prinicpal.

\section{\large \textbf{Computations and Examples}}\label{CE}

We provide several examples that justify our proposed generalizations of the 
Sawollek and other Alexander-type polynomials. In each of the examples below
we use the graded lexicographical ordering on $\tilde{\Lambda}$. Our custom
\texttt{python} code is available at \texttt{www.esotericka.org}.

\begin{example}\textup{
In~\cite{KR}, the authors give two examples of knots that are not detected 
by the generalized Alexander polynomial. One knot is the Kishino knot. The 
other is the following Kishino-like knot $K$: }
\[\raisebox{-0.5in}{
\includegraphics[height=1in]{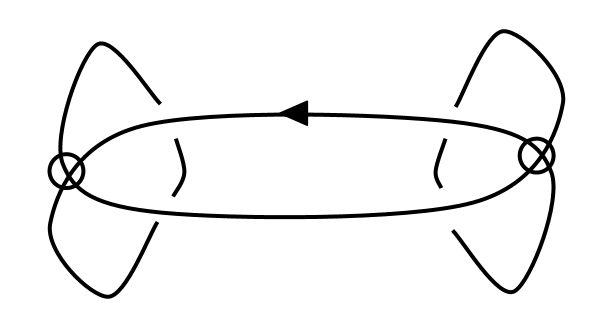}}
\]
\noindent \textup{The (standard) Kishino knot has trivial values for 
$\Delta_0^p$, $\Delta^<_0$ and $\Delta^<_1$, but this modified Kishino-like 
knot, $K$, has the following non-trivial value of $\Delta^<_1$:
}\end{example}
\begin{align*}
\Delta^<_1(K) = & \;  \{
1-t^{-1}+t^{-2}, \; \;
-1+t^{-1}+t, \; \;
-1+s, \; \; -1+s^{-1}\}
\end{align*}

\begin{example}\textup{
Two more knots that are of interest to us are Slavik's knot and 
Miyazawa's knot. Slavik's knot is not detected by the arrow 
polynomial and Miyazawa's knot is not detected by the 
Miyazawa polynomial~\cite{KD}. While Slavik's knot isn't detected by $\Delta_0^p$, 
it is detected by $\Delta^<_1$. On the other hand, $\Delta^<_1$ is 
trivial for Miyazawa's knot, but $\Delta_0^p$ and $\Delta^<_0$ are 
both non-trivial.}

\[\begin{array}{cc}
\includegraphics{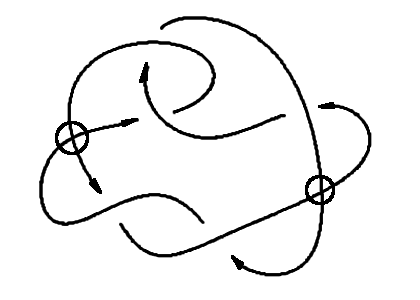} & \qquad \includegraphics{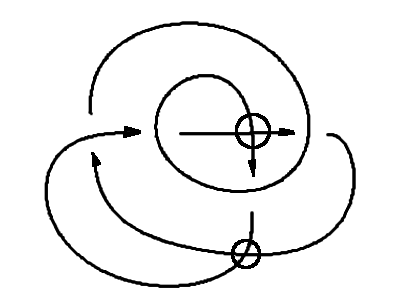} \\
\mathrm{Slavik} & \qquad \mathrm{Miyazawa} \\
\end{array}
\]
\end{example}

\noindent
\begin{align*}
\Delta^<_1(\mathrm{Slavik}) = \; &\{
3t^{-1}s - s^2 - 2t^{-2} + s^{-1}t^{-3}, \; \;
 3 - st - 2s^{-1}t^{-1} + s^{-2}t^{-2},\\
& 3s^{-1}t - t^2 - 2s^{-2} + t^{-1}s^{-3},\; \;
 -3s^{-1}t^2 + 2ts^{-2} + t^3 - s^{-3},\\\; \;
&  -3t + 2s^{-1} + st^{-2} - t^{-1}s^{-2}, \; \;
  -3s + 2t^{-1} + ts^{-2} - s^{-1}t^{-2}, \; \; \\
& -3t^{-1}s^2 + 2st^{-2}  + s^3 - t^{-3},\; \;
 -1 + tt^{-1},\; \;
 -1 + ss^{-1}\}
\end{align*}

\noindent
\begin{align*}
\Delta_0^p(\mathrm{Miyazawa}) = \; &(st-1)(s-1)(t-1)\\
\Delta^<_0(\mathrm{Miyazawa}) = \; &\{ 
-(1-s^{-1})(1-t)(s^{-1}-t), \; \;
(1 - t^{-1})(1 - s)(s - t^{-1}), \\
& (1 - t^{-1})(1 + t^{-1} - s - s^{-1}t^{-1}), \; \;
(1 - s^{-1})(1 + s^{-1} - t - s^{-1}t^{-1}) \\
& s^{-1} + t^{-1} - s - t + st - s^{-1}t^{-1},\; \;
 -1 + tt^{-1},\; \;  -1 + ss^{-1}\}.
\end{align*}

\begin{example}\textup{
Our next example considers connected sums of virtual knot diagrams. 
For classical knots, the connected sum operation is well-defined. For 
virtual knots, however, the knot type of the diagram obtained by taking 
a connected sum depends on \emph{where} the connected sum is taken. We 
give two connected sums, $K^{\#}_1$ and $K^{\#}_2$, of the virtual trefoil 
knot, 2.1, with itself:}

\[\scalebox{0.85}{$\begin{array}{ccc}
\includegraphics[height=1.1in]{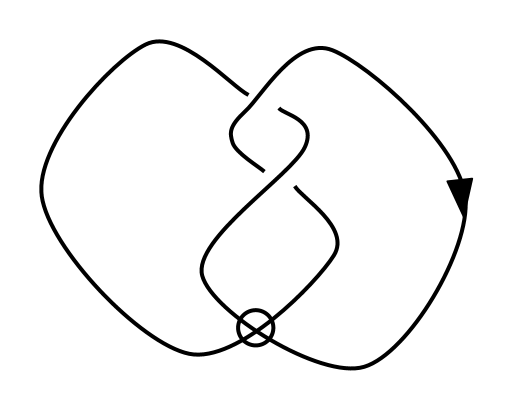}& 
\includegraphics[height=1.1in]{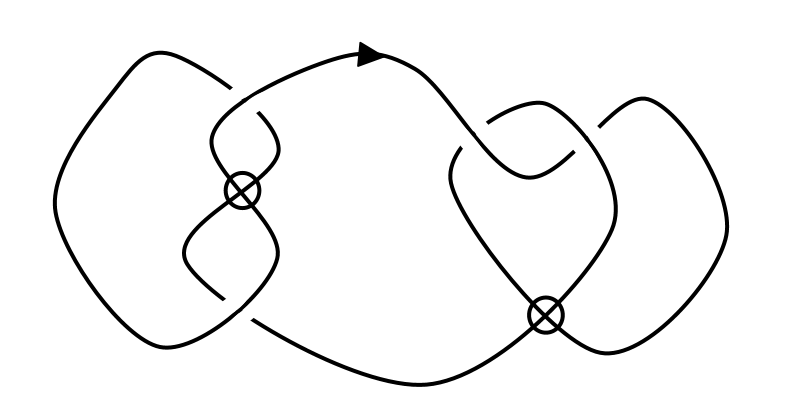} &
\includegraphics[height=1.1in]{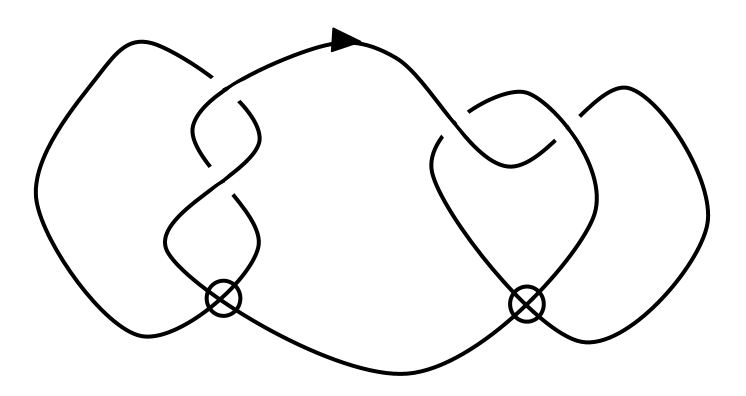}\\
\end{array}$}
\]

\noindent \textup{We see from the sets below that the invariants $\Delta_0^p$ 
and $\Delta^<_1$ distinguish these two connected sums. Note that $\Delta^<_1$ 
is trivial on both the virtual trefoil itself and on $K^{\#}_1$.
}\end{example}

\begin{align*}
\Delta_0^p(2.1) = & \; (1-s)(1-t)(1-st)\\
\Delta_0^p(K^{\#}_1) = & \; (1-s)(1-t)(1-st)(1-t+st^2+s^2t^2)\\
\Delta_0^p(K^{\#}_2) = & \; (1-s)(1-t)(1-st)(1+s-t + st^2 + s^2t^2 -ts^2 )
\end{align*}
\begin{align*}
\Delta^<_1(2.1) = & \; \{1\}\\
\Delta^<_1(K^{\#}_1) = & \; \{1\}\\
\Delta^<_1(K^{\#}_2) = & \; \{
1-t^{-1}+t^{-2}, \; \;
-1+t^{-1}+t, \; \;
-1+s, \; \; -1+s^{-1}\}
\end{align*}

For our final example, we consider the case of \textit{based virtual 
knots} or equivalently \textit{long virtual knots}; these are oriented 
virtual knots with a base point which cannot move through a classical 
crossing. In terms of diagrams, this means that strands may not move 
classically over or under or virtually detour past the base point.
Interpreting the base point as the point at infinity yields the 
long knot interpretation.

\[\includegraphics{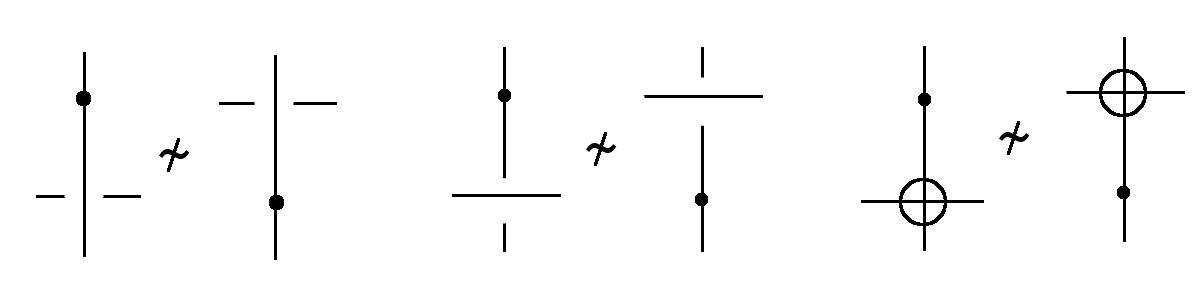}\]

Our motivation for considering based virtual knots is the observation
that while the connected sum operation is not well-defined for virtual 
knots, it \textit{is} well-defined for based virtual knots 
provided we only permit the summing operation at the base point. Indeed, based 
virtual knots form a noncommutative monoid under based connected sum, since
switching the order of summands requires forbidden moves.

The virtual biquandle of a based virtual knot with $n$ crossings has $2n+1$ generators since \emph{two} distinct generators are assigned to the semiarc containing the base point while each remaining semiarc (of which there are $2n$) is assigned a single generator.
Meanwhile, the based virtual biquandle has $2n$ relations: two coming from each of the $n$ crossings, as usual. Because the presentation matrix is no longer square, its 0-level ideal $I_0$ need not be prinicipal.
Therefore, it not surprising that the Gr\"obner invariants are stronger than
the principal invariants in the based case, as our final example shows.

\begin{example}
\textup{The two pictured oriented based virtual trefoils both have the 
trivial value $\Delta_{0}^p=1$, but are distinguished by
$\Delta_{0}^<$:}
\[\includegraphics{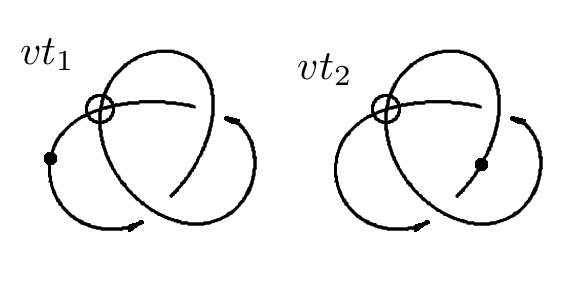}\]
\[\begin{array}{rcl}
\Delta_0^p(vt_{1}) & = & 1 \\
\Delta_0^p(vt_{2}) & = & 1 \\
\Delta_0^<(vt_{1}) & = & \{1\}\\
\Delta_0^<(vt_{2}) & = &\{1 - t^{-1} + t^{-2}, \; \;
-1 + t^{-1} + t,\; \; -1 + s,\; \; -1 + s^{-1}
\}\end{array}
\]
\textup{Taking based connected sums, we have }
\[\begin{array}{rcl}
\Delta_0^p(vt_{2}\#vt_{1}) & = &\Delta_0^<(vt_2)\\
\Delta_0^p(vt_{2}\#vt_{2}) &= & \{ -2+t+3t^{-1}-2t^{-2}+t^{-3}, \; \;
 -(1 - s^{-1})(1 - t^{-1} + t^{-2}), \\
& & 3-2t^{-1}-2t+t^{-2}+t^2,\; \;
 -1 + t^{-1}t, \\
& & (1-s^{-1})(1-t^{-1}-t), \; \;
 (1-s^{-1})^2, \; \;
\; \; s^{-1}-2+s \}\\
\end{array}\]
\end{example}

\section{\large \textbf{Questions}}\label{Q}

Many interesting questions and promising avenues of exploration involving
the $\Delta^p_k$ and $\Delta_k^<$ invariants await further work. These 
include but are not limited to:

\begin{itemize}
\item \textit{Twisted $\Delta_k$ invariants.} In the classical 
case, a matrix representation of $D_p$ coupled with a $p$-coloring of 
a knot diagram $K$ lets us replace $t$ in the Alexander 
matrix with a matrix depending on the $p$-colors
at each crossing; the principal ideal invariants of the resulting matrix 
are known as the \textit{twisted Alexander polynomials.} Matrix 
representations of labeling biquandles should be usable in an analogous way
with $AB(K)$ to define twisted $\Delta_k^p$ invariants.
\item \textit{Multivariable $\Delta_k$ invariants.} Another variant
of the Alexander polynomial for links involves multiple $t$ variables; what 
is the analogous construction for $\Delta_k$ invariants?
\item \textit{Categorification of $\Delta_k^p$.} We can apply a Khovanov-style 
construction to the state-sum expansion of $\Delta_0^p$ considered as
a determinant. How does the result compare to Knot Floer homology?
\item \textit{Skein relations.} The classical Alexander polynomial 
satisfies the well-known Conway skein relation, and the Sawollek
polynomial satisfies a similar skein relation. What skein relations,
if any, are satisfied by $\Delta_k^p$? 
\item \textit{Connected sum behavior.} We return to our original 
question: what conditions are necessary and sufficient for a connected
sum of virtual knot or link diagrams to behave like classical knots 
under connected sum with respect to the $\Delta_k^p$ invariants?
A related question: what is the center of the monoid of based virtual
knots?
\item \textit{Virtual links with boundary.} Based virtual knots are a special
case of virtual knots with boundary, where the supporting surface of the
virtual knot has a fixed boundary and knots or links on the surface may have
endpoints in the boundary. Gluing such surfaces along boundary components
such that endpoints of knots match up generalizes our based connected sum 
operation. What is the algebraic structure of such knots? The authors are 
grateful to Charlie Frohman for this observation.
\end{itemize}


%
%
%
%

\end{document}